\documentclass[11pt]{amsart}
\usepackage{amsmath}
\usepackage{amssymb}
\usepackage{amscd}
\newtheorem{thm}{\sc theorem}[section]
\newtheorem{dfn}[thm]{\sc definition}

\newtheorem{lem}[thm]{\sc lemma}
\newtheorem{cor}[thm]{\sc corollary}

\newcommand\prf[1]{\par\vspace*{-.5ex}\noindent{\it Proof {#1} --- }}
\def\sur#1#2{{\buildrel{#1}\over{#2}{}}}
\def\germ#1#2#3{{{#1}\sur{#2}\leftarrow{#3}}}
\newcommand{\cqfd}{\hfill{$\bullet$}\par\vspace{3ex}}

\newcommand{\F}{\mbox{$\mathcal F$}}

\newcommand{\rank}{\mbox{$\mbox{\rm rank}$}}

\def\N{{\bf N}}
\def\Z{{\bf Z}}
\def\Q{{\bf Q}}
\def\R{{\bf R}}

\def\S{{\bf S}}

\def\D{{\bf D}}

\def\mun{{^{-1}}}

\begin{document}

\vskip 1cm
\today

\title{a compactly generated pseudogroup which is not realizable}

\author{\rm Ga\"el Meigniez}
\footnote{Laboratoire de Math\'ematiques et d'Application des Math\'ematiques,
Universit\'e de Bretagne Sud,
Universit\'e Europ\'eenne de Bretagne. Supported by the Japan Society for Promoting Science,
No.\ L03514.}
\address{
Universit\'e de Bretagne-Sud, L.M.A.M.,
BP 573, F-56017 Vannes, France}
\email{Gael.Meigniez@univ-ubs.fr}



\thanks{This work was made while the author was visiting the University of Tokyo. It is a pleasure to thank specially T. Tsuboi, and the animators of the \emph{Saturday Seminar,} for a hospitality particularly favorable to research.}

\begin{abstract} We exhibit a pseudogroup of smooth local transformations of the real line, which is compactly generated, but not realizable as the holonomy pseudogroup of a foliation of codimension 1 on a compact manifold. The proof relies on a description of all foliations with the same dynamic as the Reeb component.
\end{abstract}

\maketitle

\thispagestyle{empty}
\vskip 1cm
\section{Introduction}

To every foliated manifold $(M,\F)$
of arbitrary dimension and codimension, one associates, following Ehresmann, a {\it pseudogroup} $Hol(\F)$ of local
transformations, called the {\it holonomy pseudogroup,}
that represents
its ``dynamic'' or ``transverse structure''. The holonomy pseudogroup is well-defined up to some natural equivalence
between pseudogroups: {\it Haefliger equivalence.}

The inverse {\it realization} problem has been raised: make a foliation with prescribed dynamic,
 the ambiant manifold not being prescribed ( but it must be compact.)
More precisely, given a pseudogroup $G$, make if possible
a {\it compact} foliated manifold $(M,\F)$
such that $Hol(\F)$ is Haefliger-equivalent to $G$~.

Note that if one drops the compactness condition, the question vanishes~:
every pseudogroup is easily seen to represent the dynamic of some foliated {\it open} manifold.

Andr\'e Haefliger made the realization problem precise by exhibiting a necessary condition, {\it compact generation.}
Is it sufficient?
Some partial positive answers have been given for rather rigid
species of pseudogroups \cite{cav2002}\cite{hae1984}\cite{mei1995}\cite{mei1997}.

The object of this paper is to answer negatively in general.
 We give a counterexample among pseudogroups of smooth local transformations of the real line.

\section{Pseudogroups}
We recall briefly the definitions and the basic properties, now classical.
See also \cite{hae2002}.
\medbreak
 An arbitrary differentiability class is understood. Let $T$ be a manifold, not necessarily compact. A boundary is allowed.

A \emph{local transformation} of $T$ is a diffeomorphism between two nonempty open subsets $Dom(\gamma)$~,
$Im(\gamma)$ of $T$~. The compose $\gamma'\gamma$ is defined whenever $Im(\gamma)$ meets $Dom(\gamma')$~, and its
domain is $\gamma\mun(Dom(\gamma'))$~. If $A$~, $B$ denote two sets of local transformations, then as usual $AB$
denotes the set of their composes. If $U\subset T$~, then $A\vert U$ denotes the set of the elements of $A$ whose domains
and images are both contained in $U$~.

Note that $\partial T$ is necessarily preserved by every local transformation.

\begin{dfn}\cite{veb1931}
A \emph{pseudogroup} on $T$ is a set $G$ of local transformations such that~:
\begin{enumerate}
\item For every nonempty open $U\subset T$ the identity map $1_U$ belongs to $G$~;
\item $GG=G\mun=G$~;
\item For every local transformation $\gamma$ of $T$~, if $Dom(\gamma)$ admits an open cover $(U_i)$
 such that every restriction $\gamma\vert U_i$
belongs to $G$~, then $\gamma$ belongs to $G$~.
\end{enumerate}
\end{dfn}

Then, by (1) and (2), $G$ is also stable by restrictions: if $\gamma$ belongs to $G$ and if $U\subset Dom(\gamma)$ is nonempty open, then $\gamma\vert U$ belongs to $G$~.

For example, every set $S$ of local transformations of $T$ is contained in a smallest pseudogroup $<S>$ containing $S$~. Call $<S>$ the pseudogroup
\emph{generated} by $S$~.
\medbreak
Every point $t$ in $T$ has under a pseudogroup $G$~:
\begin{enumerate}
\item
An \emph{orbit}: the set of the images $\gamma(t)$ through
 the local transformations
 $\gamma\in G$ defined at $t$~;
\item
An \emph{isotropy group}: the group of the germs at $t$ of the
 local transformations $\gamma\in G$ defined at $t$ and fixing $t$~.
\end{enumerate}
\medbreak
Let $(M,\F)$ be a manifold foliated in codimension $q$~. By a \emph{transversal} one means a $q$-manifold $T$ immersed into $M$
transversely to $\F$~, not necessarily compact, and such that $\partial T=T\cap\partial M$~. One calls $T$ \emph{exhaustive} (or \emph{total}) if it meets every leaf.

\begin{dfn}\cite{ehr1954}
The holonomy pseudogroup $Hol(\F,T)$ of a foliation $\F$ on an exhaustive transversal
 $T$ is the pseudogroup \emph{generated} by the local transformations of $T$ of the form $f(x,0)\mapsto f(x,1)$~,
where~:
\begin{enumerate}
\item $D^q$ is the open $q$-disk~;
\item $f:D^q\times[0,1]\to M$ is a map transverse to $\F$~;
\item $f^*\F$ is the foliation on $D^q\times[0,1]$ by the first projection;
\item $f$ embeds $D^q\times 0$ and
$D^q\times 1$ into $T$~.
\end{enumerate}
\end{dfn}

This holonomy pseudogroup does represent the dynamic of the foliation in the sense that
 there is a one-to-one correspondance $L\mapsto L\cap T$ between the leaves of $\F$ and  the orbits of $Hol(\F,T)$~.
 A topologically closed orbit
corresponds to a closed leaf. The isotropy group of $Hol(\F,T)$
at any point
  is isomorphic with the holonomy group of the corresponding leaf. Etc.
\begin{dfn}\cite{hae1984}
A \emph{Haefliger equivalence} between two pseudogroups $(T_i,G_i)$ ($i=0,1$) is a pseudogroup $G$
on the disjoint union of $T_0$ with $T_1$~, such that
 $G\vert T_i=G_i$ ($i=0,1$) and that no orbit of $G$ in entirely contained
 in $T_0$ or in $T_1$~.
\end{dfn}

For example, obviously, the two holonomy pseudogroups of a same foliation on two exhaustive transversals are Haefliger equivalent.

A Haefliger equivalence between $(T_1,G_1)$ and $(T_2,G_2)$
 induces a one-to-one correspondance between the orbit spaces
 $T_i/G_i$ ($i=0,1)$~;
 a closed orbit corresponds to a closed orbit; the isotropy groups at points on corresponding orbits are isomorphic; etc.
\medbreak
Let $(T,G)$ be a pseudogroup. Call a subset $T'\subset T$ \emph{exhaustive} if it meets every orbit.
Call $\gamma\in G$ \emph{extendable} if it is the restriction to $Dom(\gamma)$ of some
 $\bar\gamma\in G$ such
 that $Dom(\gamma)$ is relatively
compact in $Dom(\bar\gamma)$~.

\begin{dfn}\cite{hae1985} A pseudogroup $(T,G)$ is \emph{compactly generated} if
 there are an exhaustive, relatively compact, open subset $T'\subset T$~, and finitely many elements
of $G\vert T'$ which are {extendable} in $G$ and which generate $G\vert T'$~.
\end{dfn}

 This property is invariant by Haefliger equivalence \cite{hae1985}\cite{hae2002}.
The holonomy pseudogroup of every foliated compact manifold is compactly generated \cite{hae1985}. Also, recently N. Raimbaud has given
 a natural generalization of compact generation in the realm of Lie groupoids, where it is a Morita-equivalence invariant \cite{rai2009}.

\section{The example}
In this paper, to fix ideas one works in the smooth ($C^\infty$)
 differentiability class; all foliations and pseudogroups are
 transversely orientable; all diffeomorphisms are orientation-preserving.

In the realization problem, one may allow that $M$ have some boundary components transverse to $\F$~, or not. This has no influence on the answer. Indeed assume that some pseudogroup $G$ is realized by $(M,\F)$ which has some
 transverse boundary components $\partial_{tr}M$~. Let $\D^2$ denote
 the compact 2-disk. Then $G$ is also
 realized by $(M',pr_1^*\F)$ where $M'$ is, in $\partial(M\times\D^2)$~,
 the union of $M\times\S^1$ with $\partial_{tr}M\times\D^2$~.

\medbreak
 The counterexample to realizability is as follows.
Let $\alpha$~, $\beta$ be two global diffeomorphisms of the real line such that~:
\begin{enumerate}
\item $\alpha$ is
 a contraction fixing $0$~,
 that is, $\vert\alpha(t)\vert<\vert t\vert$ for every $t\neq 0$~;
\item The support of $\beta$ is compact and contained in $(-1,0]$~;
\item The germs
of $\alpha$ and $\beta$ at $0$ generate a nonabelian free group.
\end{enumerate}
The third condition is in some sense generically fulfilled;
one can also make an explicit example by the following classical method.

Let $A, B\in PSL(2,\R)$ generate a free group.
Lift them into two diffeomorphisms $\tilde A$~, $\tilde B$ of the real line commuting with the unit translation.
 Composing if necessary $\tilde A$ with some integral translation, $\tilde A(t)>t$ for every $t$~. After a conjugation
 by the exponential map, one has two diffeomorphisms $a$~, $b$ of $(0,+\infty)$ generating a free group.
Moreover they verify the tameness property:$$Ct\le a(t),b(t)\le C't$$for some constants $0<C<C'$~.
After a new
conjugation by $\phi:t\mapsto\exp(-\exp(1/t))$~,
one has two germs of diffeomorphisms $f:=\phi\mun a\phi$ and $g:=\phi\mun b\phi$ on the right-hand side of $0$~. It is easily verified that~:
$$\vert{\phi\mun(Ct)\over\phi\mun(t)}-1\vert=o(\phi\mun(t)^n)$$
for every $C>0$~, $n$~, and $t\to 0$~. Thus $f$ and $g$ are
 flat on the identity at 0; and it remains only to change the orientation on the line, and to extend both germs in an obvious way, to get to diffeomorphisms $\alpha$~, $\beta$ with the prescribed properties.

\begin{thm}\label{counterexample_thm}The pseudogroup
 $G:=<\alpha,\beta>$ generated by the above
 diffeomorphisms
 is compactly generated and is not realizable~.
\end{thm}

The first affirmation is actually easy:

\begin{lem}{ $G$
 is compactly generated.}\end{lem}

{\it Proof ---}
Take $T':=(-1,1)$ and $\alpha':=\alpha\vert T'$ and $\beta':=\beta\vert T'$~.
 Then obviously $T'$ is exhaustive
and $\alpha'$~, $\beta'$ are extendable in $G$~. It remains to verify that
they do generate $G\vert T'$~.

Let be given the germ, denoted
 $\germ{\gamma(t)}{\gamma}t$~, of some element $\gamma$ of $G\vert T'$ at some
point $t$ in its domain. Thus $t$~, $\gamma(t)\in T'$~. We have to write
this germ as a compose of germs of the diffeomorphisms $\alpha^{\pm 1}$
and $\beta^{\pm 1}$ \emph{all taken at points of $T'$~} --- and this is the marrow of bone of compact generation.

 But here it is easy: since $G$ is
 generated by $\alpha$ and $\beta$~, by definition the given germ decomposes as a
 composable sequence~:
$$(\germ{\gamma(t)}\gamma t)=\big(\germ{t_n}{\gamma_n}{}\germ{t_{n-1}}{\gamma_{n-1}}{\dots}\germ{t_1}{\gamma_1}{t_0}\big)$$
 of germs of $\alpha^{\pm 1}$
and $\beta^{\pm 1}$ at some points $t_0=t$~, \dots, ${t_n}\in\R$~.

 Take such a decomposition of minimal length $n$~. Then we claim that  $t_0$~,\dots, ${t_n}\in T'$~. Indeed,
if not, one has for example$$t_\ell:=\sup\{t_0,\dots,t_n\}\ge 1$$

By maximality of $t_\ell$~, and since $\alpha(t_\ell)<t_\ell$~,
 one has either $\gamma_\ell=\beta^{\pm 1}$ or $\gamma_{\ell+1}=\beta^{\pm 1}$
or $\gamma_\ell\mun=\alpha=\gamma_{\ell+1}$~, contrarily to the minimality
of the length of the decomposition.\cqfd

\def\RR{{\mathcal R}}

Observe that the halfline $[0,+\infty)$ is saturated for $G$~, and that the restriction $G\vert[0,+\infty)$
is actually the transverse structure of a Reeb component. The proof that $G$ is not
 realizable will rely on a precise description of all the
foliations with the same transverse structure as a Reeb component,
 from which it will then follow that
the boundary leaf cannot present such
a free holonomy group on the side exterior to the component.
\section{Generalized Reeb components}

Fix a contraction $\eta$ of the halfline
$\R_+:=[0,+\infty)$~, and consider the generated pseudogroup $<\eta>$~.

Of course, this pseudogroup has a canonical
realization in dimension $3$~: the classical Reeb foliation
on $\D^{2}\times\S^1$~, obtained as follows.
Having foliated $(\R^{2}\times\R_+)\setminus 0$ by
 its projection
 onto the halfline,
one passes to the quotient by
the foliation-preserving diffeomorphism $(x,t)\mapsto
(x/2,\eta(t))$~.

This obvious construction has a natural generalization
(Alca\-lde-\-Cuesta-\-He\-ctor-\-Schwei\-tzer, unpublished). One is
given a compact connected $(n-1)$-manifold
 $C$ with smooth connected boundary and
a self-embedding$$\phi:C\to Int(C)$$(In the classical case,
 $C$ would be $\D^{2}$~; and $\phi(x)=x/2$~).
 From these data, one makes a \emph{generalized Reeb component} as follows.

Consider the projective limit~:
\[P:=\cap_{i\in\N}\phi^i(C)\]
and the inductive
limit~:
\[I:=(C\times\Z)/((x,i+1)\sim(\phi(x),i))\]
Denote $[x,i]$ the class of the pair $(x,i)$~. One has
a diffeomorphism~:
\[\Phi:I\to I:[x,i]\mapsto[\phi(x),i]\]
Identify $C$ with a subset of $I$ through the embedding $x\mapsto[x,0]$~. Thus $\Phi\vert C=\phi$~. It is also convenient
to fix a smooth function $f_0$ on $C\setminus Int(\phi(C))$ such that $f_0\mun(0)=\partial C$ and that $f_0\mun(1)=\phi(\partial C)$~. It extends uniquely into a function $f$ on $I\setminus P$ such that $f\circ\Phi=f+1$~. Obviously, $f$ is proper.
 Set $f=+\infty$ on $P$~.

Also, let $g$ be a function on $(0,+\infty)$ such that $g(\eta(t))=g(t)+1$~. Set $g(0)=+\infty$~.

Define~:
\[\tilde R:=(I\times\R_+)\setminus(P\times 0)\]
Foliate it by its projection onto $\R_+$~. Also endow it with the foliation-preserving
 diffeomorphism~:\[\gamma:\tilde R\to\tilde R:(x,t)\mapsto
(\Phi(x),\eta(t))\] and with the function~:
\[F(x,t):=\min\{f(x),g(t)\}\]
It is immediately verified that $F$ is finite and proper, and that $F\circ\gamma=F+1$~.
It follows that $\gamma$
acts freely, properly discontinuously and cocompactly on $\tilde R$~.
Thus the quotient is a foliated, compact, Hausdorff manifold  $(R,{\mathcal R})$~.
\begin{dfn}\label{DGRC}
Call $(R,{\mathcal R})$ the generalized Reeb component associated to
the self-embedding $(C,\phi)$~.
\end{dfn}

Obviously  $(R,{\mathcal R})$ realizes $<\eta>$~. Conversely~:
\begin{thm}\label{GRC}
Every realization of the pseudogroup generated by a  contraction of the halfline
is diffeomorphic to some generalized Reeb component in the sense of Alcalde-Cuesta-Hector-Schweitzer.
\end{thm}

Here ``realization'' is understood without transverse boundary components.
\medbreak
\prf{of theorem \ref{GRC}} Given a realization  $(M,\F)$ of the contraction $\eta$~,
\begin{lem}\label{dev} The foliation $\F$
is developable over $\R_+$ and complete~\dots
\end{lem}\dots which means the following~:
there is an infinite cyclic covering~:\[\pi:\tilde M\to M\]
a generator $\gamma$ of the deck transformation group,
and a ``developing map''~:
\[D:\tilde M\to\R_+\]
such that~:
\begin{enumerate}
\item The map $D$ is a surjective submersion;

\item The fibres of $D$ are connected and are the leaves of $\pi^*\F$~;

\item One has
$\eta\circ D=D\circ\gamma$~.
\end{enumerate}
\prf{of lemma \ref{dev}} One can either call to the general theory of transversely affine foliations, or deduce these properties from the corresponding ones observed on an explicit classifying space, as follows.

Changing $3$ and $2$ into $+\infty$ in the above
 construction of the classical Reeb component,
one gets an infinite-dimensional Reeb component
$(R^\infty,\RR)$~. The holonomy
covering of each leaf is weakly contractible. That foliation
thus being the classifying space of its pseudogroup $<\eta>$ \cite{hae1984}, there exists a
 classifying map $c:M\to R^\infty$ transverse to  $\RR$
 such that $\F=c^*\RR$~, and that
$c$ induces a Haefliger equivalence between the holonomy pseudogroups
 of $\F$ and
 of $\RR$~.

In particular $c$ induces a bijection of the leaf spaces~; and,
for every leaf $L$ of $\RR$~, the map $c$
also induces a group isomorphism from the
holonomy group of the leaf $c\mun(L)$ onto  the
holonomy group of $L$~.

Thence $c$ maps the holonomy group of $\partial M$
onto the holonomy group of $\partial R^\infty$~. Thus $c$ maps the fundamental
group $\pi_1M$ onto $\pi_1R^\infty\cong\Z$~, hence an infinite cyclic
covering $\tilde M$ and a lifting $\tilde c:
\tilde M\to\tilde R^\infty$~.
Define $D$ as $\tilde c$ followed by the projection to $\R_+$~
. The above properties of $c$ immediately translate into the demanded properties for $D$~.
\cqfd

\def\NN{{\mathcal N}}
(Continuation of the proof of theorem \ref{GRC}) Fix
 in $M$ an arbitrary smooth foliation $\NN$ of dimension 1 transverse to $\F$~. In particular $\NN$ is
 transverse to $\partial M$~. Lift it into
a foliation  $\tilde\NN$ of the covering $\tilde M$~.
Consider the canonical projection onto the space of orbits:
\[pr:\tilde M\to I:=\tilde M/\tilde\NN\]and
 the homeomorphism $\Phi:I\to I$ such that:
\[pr\circ\gamma=\Phi\circ pr\]
and the $\Phi$-invariant,
topologically closed subset:\[P:=I\setminus pr(\partial\tilde M)\]
\begin{lem}\label{product} The space of orbits $I$ is a connected Hausdorff manifold.
Moreover, there is a diffeomorphism~:
\[\tilde M\cong(I\times\R_+)\setminus(P\times 0)\]
through which $\gamma(x,t)=(\Phi(x),\eta(t))$ and $D(x,t)=t$ and $pr(x,t)=x$~.
\end{lem}
\prf{}
The halfline bears an $\eta$-invariant vector field $u(t)\partial/\partial t$~, smooth and nonsingular in $(0,+\infty)$~, null at 0.
It needs not be differentiable at 0. Clearly it is complete.  Let $(\eta^s)_{s\in\R}$ be the associated 1-parameter
group of homeomorphisms of the halfline.
Consider the unique vector field $\tilde X$
 in $\tilde M$ tangent to $\tilde\NN$ and
projecting onto $u(t)\partial/\partial t$ through $D$~. Since $\tilde\NN$ is $\gamma$-invariant, since $u(t)\partial/\partial t$ is
  $\eta$-invariant and since $D$ is equivariant, $\tilde X$ is $\gamma$-invariant. In other words $\tilde X$ is the pullback into $\tilde M$ of some
 vector field $X$ on the compact manifold $M$~, which is smooth in the interior of $M$ and null on $\partial M$~. The vector field
 $u(t)\partial/\partial t$ being complete, one concludes easily that $X$ is complete. Thus
$\tilde X$ is complete. Let $(\xi^s)_{s\in\R}$ be the associated 1-parameter
group of homeomorphisms of $\tilde M$~.
 Then $D\circ\xi^s=\eta^s\circ D$~. From this equivariance
follows easily that the following map is one-to-one and onto~:
\[\psi:D\mun(1)\times\R_+^*\to Int(\tilde M):(x,t)\mapsto\xi^{\int_1^t{d\tau\over u(\tau)}}(x)\]
Being obviously etale, it is a diffeomorphism.

In particular $I$
is diffeomorphic to $D\mun(1)$~, thus a connected Hausdorff manifold.

It remains to extend $\psi$ to the boundary. For every $x\in
I\setminus P$~, set
$\psi(x,0):=\lim_{s\to-\infty}\xi^s(x)\in \partial\tilde M$~.
Obviously this extends $\psi$ into a global diffeomorphism from $(I\times\R_+)\setminus(P\times 0)$ onto
$\tilde M$~,
through which $\gamma(x,t)=(\Phi(x),\eta(t))$ and $D(x,t)=t$ and $pr(x,t)=x$~.
\cqfd

It seems that a little more work is necessary to make the dynamic
of $\Phi$~, and its relation to $P$~, precise; and thus to achieve the proof
of theorem \ref{GRC}. For example, at this point it is not obvious that $P$ is compact.

 One identifies $\tilde M$ with $(I\times\R_+)\setminus(P\times 0)$~.

 It is a well-known property of infinite cyclic coverings that $\tilde M$
admits a proper smooth function $F$ such that $F\circ\gamma=F+1$~. To fix ideas, one can arrange that $0$
is a regular value of $F$ and of $F\vert\partial\tilde M$~. Also, by \ref{dev},
 $\partial\tilde M$ is connected. Thence one can arrange also that $\partial(F\mun(0))$ is connected.

\begin{lem}\label{uslim} (i) For every $x\in I$~, one has~:
\[\lim_{t\to+\infty}F(x,t)=-\infty\]

(ii)  For every $p\in P$~, one has~:
\[\lim_{t\to 0}F(p,t)=+\infty\]
(iii) More precisely, for every $p\in P$~, one has~:
\[\lim_{(x,t)\in\tilde M,(x,t)\to (p,0)}F(x,t)=+\infty\]
\end{lem}
\prf{}
 (i) Let $T$ be the maximum of $D(x,t)=t$ on the compact fundamental domain
 $F\mun([-1,0])$~. Let $g$ be a decreasing function on $(0,+\infty)$ such that $g(T)=0$ and $g\circ\eta=g+1$~. Then $F(x,t)\le g(t)$ at every point $(x,t)$ of $\tilde M$~. Thus $F(x,t)\to-\infty$ for $t\to+\infty$~.

(ii) The halfline $p\times(0,1]$ being properly embedded in
$\tilde M$~, the limit exists, either $-\infty$
or $+\infty$~.
By contradiction, assume that it is $-\infty$~.
For every $i$ large enough~:
\[F(\gamma^i(p,1))=F(p,1)+i>0\]
thus the halfline $\gamma^i(p\times(0,1])=\Phi^i(p)\times(0,\eta^i(1)]$ would meet $F\mun(0)$ in at least one point
 $(\Phi^i(p),t_i)$~. The level set $F\mun(0)$ being compact, some subsequence of the sequence $(\Phi^i(p),t_i)$ converges to
some $(q,t)\in F\mun(0)$~.
Since $t_i\le\eta^i(1)$~,
one has $t=0$~. Since  $P$ is $\Phi$-invariant and topologically closed in $I$~, one has $q\in P$~. Thus $(q,t)\in P\times 0$~,
the desired contradiction.

(iii) Consider a fundamental sequence $(V_i)$ of connected neighborhoods of $p$ in $I$~, and~:
\[W_i:=(V_i\times[0,1/i])
\cap\tilde M\]
and fix a large positive $T$~. Since $F\mun[-T,+T]$ is compact and does not contain $(p,0)$~, it is disjoint from $W_i$
for every $i$ large enough. Since $W_i$ is connected, either $F>T$ on $W_i$ or $F<-T$ on $W_i$~. The second possibility
being ruled out by (ii), the lemma is proved.
\cqfd

On $I\setminus P$ one has the proper function $f(x):=F(x,0)$ and one defines~:\[C:=P\cup f\mun[0,+\infty)\subset I\]
\begin{cor} The subset
$C\subset I$ is a compact submanifold of codimension $0$ with smooth boundary and $P$ is contained in its interior.
Both $C$ and $\partial C$ are connected.
\end{cor}
\prf{}
By lemma \ref{uslim}, firstly {\it $C$ is relatively compact in $I$~.} Indeed, for every $x\in C$~,
by (i) and (ii) the halfline $x\times\R_+$ meets $F\mun(0)$~. That is, $C$ is contained in $pr(F\mun(0))$ which
is compact in $I$~.

Secondly, $P$ is contained in the topological interior of $C$~. This follows at once from (iii).
In particular, the topological boundary of $C$ in $I$ is $f\mun(0)=\partial F\mun(0)$~,
 a smooth compact connected $(n-2)$-manifold.
Since $I$ and $\partial C$ are connected, $C$ is connected.
 \cqfd
Now, recalling that one has the diffeomorphism $\Phi$ of $I$ such that $\Phi(P)=P$ and that
 $f\circ\Phi=f+1$ on $I\setminus P$~, one gets easily~: 
\[\Phi(C)\subset Int(C)\ \ \rm{and}\ \ P=\cap_{i\in\Z}\Phi^i(C)\ \ \rm{and}\ \ I=\cup_{i\in\Z}\Phi^i(C)\]
By lemma \ref{product}
 the foliated manifold
$(M,\F)$ is diffeomorphic to the generalized Reeb component
associated with $(C,\Phi\vert C)$ according to definition \ref{DGRC};
 and the theorem \ref{GRC} is proved.
\cqfd
\smallbreak
In general, the cobordism $C\setminus Int(\phi(C))$ is of course not trivial. Accordingly,
 the boundary leaf $\partial R$ of an arbitrary
generalized Reeb component (definition \ref{DGRC})
is not necessarily fibred over the circle. However, we always have the following finiteness property,
well-known e.g. in the classical study of knots~:
\begin{lem}\label{finiteness} Let $\partial\tilde R$ be the holonomy covering
of the boundary leaf of a
generalized Reeb component. Then the homology groups of
$\partial\tilde R$ with coefficients in any field $k$ are of finite rank over $k$~.
\end{lem}
\prf{} (All homology groups are with coefficients in $k$~.) Let $C$~, $\phi$~, $I$~, $P$~, $\Phi$~, $R$
 be as in definition \ref{DGRC}. For every positive $i$~,
write $C_i=\Phi^{-i}(C)$
and $W_i=C_i\setminus Int(C_{-i})$~. In the following
commutative diagram (where
all arrows are induced by inclusions)~:
\[
\begin{CD}
H_*(C_i)  @>{}>> H_*(C_i,C_i\setminus Int(W_i)) \\
@A{}AA @A{\rho}AA\\
H_*(W_i)  @>{\beta}>>  H_*(W_i,\partial W_i)
\end{CD}
\]
the right-hand vertical arrow $\rho$
 is invertible by the excision theorem, thus~:
\[\rank(\beta)\le\rank H_*(C_i)\]
On the other hand, the long exact relative homology sequence
for the couple $(W_i,\partial W_i)$ gives~:
\[\rank H_*( W_i)\le\rank(\beta)+\rank H_*(\partial W_i)\]
But $C_i$ is diffeomorphic to $C$ and $\partial W_i$ is diffeomorphic to
two copies of $\partial C$~, thus:
\[\rank H_*(W_i)\le\rank H_*(C)+2\rank H_*(\partial C)\]
An upper bound independant on $i$~.
The covering space $\partial\tilde R$ being
 the inductive limit of the sequence~:
\[W_1\subset W_2\subset\dots\subset W_i\subset\dots\]
the rank of $H_*(\partial\tilde R)$ admits the same majoration.
\cqfd

\section{Proof of theorem \ref{counterexample_thm}} Consider again
the pseudogroup $G=<\alpha,\beta>$~, where $\alpha$ is a contraction of
the real line $\R$ fixing 0 and where $\beta$ is a diffeomorphism of $\R$ with compact
support contained in $\R_-$~, and such that their germs at 0
generate a nonabelian free group. In the pseudogroup $(\R,G)$ one may call $\R_+$ a
\emph{paradoxical Reeb component:} a saturated domain with the same dynamic as a Reeb component, but whose boundary $0$
has a complicated isotropy group outside.

On the contrary, the preceding section has shown us that
the corresponding paradoxical
 Reeb components cannot
 exist among foliations,
and so $G$ is not realizable.

More precisely,
in the isotropy group $Iso(G,0)$ of $G$ at point $0$~,
 one has the subgroup $ExtIso(G,0)$ consisting of
the germs which are the identity on the right-hand side of $0$~.
 Clearly $ExtIso(G,0)$
is the normal subgroup generated by $\beta$~,
and thus a nonabelian free group of infinite rank. Consider its abelianization
(quotient by the derived subgroup) $ExtIso(G,0)_{ab}$~.
Then the vector space$$\Q\otimes ExtIso(G,0)_{ab}$$
 is of infinite rank
over $\Q$~.

On the other hand, assume by contradiction that $G$ has
 some realization $(M,\F)$~.
 That is, $(M,\F)$ would be a foliated compact manifold
 whose holonomy pseudogroup would be
 Haefliger-equivalent to $G$~.
 As aforesaid, one can assume moreover that $M$ is closed.
The halfline $\R_+$ being
 $G$-invariant, $M$ would contain a compact saturated domain $R$
that would realize the pseudogroup $G\vert\R_+$~, that is,
 the pseudogroup on the half line generated by the contraction $\alpha$~.
After
theorem \ref{GRC}, $R$ would be a generalized Reeb component.
 Let $Hol(\F,\partial R)$ denote the holonomy
group of the leaf $\partial R$~,
 and $ExtHol(\F,\partial R)$ denote the subgroup
of germs which are the identity inside $R$~.
Let also $\tilde\partial R$ be the infinite cyclic covering corresponding
to the holonomy inside $R$~.
 So, $\pi_1\tilde\partial R$ is mapped onto $ExtHol(\F,\partial R)$~.
In consequence, the vector space$$\Q\otimes ExtHol(\F,\partial R)_{ab}$$
being a quotient of $H_1(\tilde R;\Q)$~,
 which is of finite rank after lemma \ref{finiteness},
 is also of finite rank
over $\Q$~.

But, since the holonomy pseudogroup of $\F$
 is Haefliger-equivalent to $G$~, the groups $ExtHol(\F,\partial R)$
 and $ExtIso(G,0)$ are of course isomorphic,
a contradiction.\cqfd
\section{Questions}

Haefliger has introduced an interesting stronger notion of \emph{compact presentation} for pseudogroups \cite{hae2002}.
 The holonomy pseudogroup of any foliated compact manifold is compactly presented, and any compactly presented pseudogroup
is compactly generated. Unfortunately, compact presentation seems difficult to decide on explicit examples such as ours.

{\it Question ---} Is the above pseudogroup $<\alpha,\beta>$ compactly presented~?

Presently, I know no pseudogroup which is compactly generated but not compactly presented.
\medbreak
In a direction complementary to the present paper, in a forthcoming one I will show that actually many compactly generated pseudogroups of codimension 1 are realizable, and even realizable on manifolds of small dimension. The result is as follows.

 Let $(T,G)$ be a compactly generated pseudogroup, with $\dim T=1$~. The notion of ``dead end component'', well-known for codimension one foliations, has an obvious analogue for pseudogroups. Those components are bounded by closed orbits, of which we consider the isotropy groups. One can show:
\medbreak
{\it 
1. If every dead end boundary isotropy group is solvable,
 then $(T,G)$ is realizable in a 4-manifold.

2.  $(T,G)$ is realizable in a 3-manifold if and only if every dead end boundary isotropy group is
abelian of rank $\le 2$~.}
\medbreak
In particular, if $G$ has no closed orbit, or more generally no dead end component, then it is realizable in dimension 3.
If $G$ is $PL$~, or projective (local transformations of the type $t\mapsto (at+b)/(ct+d)$~)
 then it is realizable in dimension 4.

{\it Question ---} Is every real-analytic compactly generated pseudogroup of codimension 1 realizable?
\medbreak
So, one has seen in the present paper a sufficent condition for \emph{not} being realizable
(some Reeb component boundary isotropy group is nonabelian free)
and one will see also a sufficent condition for being realizable (every dead end boundary isotropy group is solvable).
These two conditions are not exactly complementary, there remains a little gap. Maybe a good understanding of
 {\it compact presentation} in codimension 1 would fill
the gap.

\end{document}